\documentstyle[11pt]{article}
\newcommand{\too}{\longrightarrow}

\newcommand{\U}{{\cal U}}
\newcommand{\V}{{\cal V}}
\newcommand{\Li}{{\cal L}}

\newcommand{\J}{{\cal J}}\newcommand{\h}{{\cal H}}

\newcommand{\C}{{\cal C}}
\newcommand{\X}{{\cal X}}
\newcommand{\G}{{\cal G}}
\newcommand{\D}{{\cal D}}

\newcommand{\di}{\displaystyle}
\newcommand{\Om}{\Omega}
\newcommand{\na}{\nabla}
\newcommand{\wi}{\widetilde}
\newcommand{\al}{\alpha}
\newcommand{\be}{\beta}
\newcommand{\ga}{\gamma}
\newcommand{\Ga}{\Gamma}
\newcommand{\e}{\epsilon}
\newcommand{\la}{\lambda}
\newcommand{\De}{\Delta}

\def\reel{ {\rm I}\!{\rm R} }

\def\prs{<\;,\;>}
\def\br{[\;,\;]}

 \def \rat{ {\rm Q}\kern-.65em {}^{{}_/ }}

\newtheorem{Def}{Definition}[section]
\newtheorem{th}{Theorem}[section]
\newtheorem{pr}{Proposition}[section]
\newtheorem{Le}{Lemma}[section]
\newtheorem{co}{Corollary}[section]
\newtheorem{rem}{Remark}[section]
\newtheorem{exem}{Example}[section]

\title{  Riemannian Geometry of Lie algebroids  }
\author{Mohamed Boucetta} \date{ }\parindent=0cm
\begin{document}
\maketitle

{\bf Abstract.} We introduce Riemannian Lie algebroids  as a
generalization of Riemannian manifolds and we show that   most of
the classical tools and results known in Riemannian geometry can
be stated in this setting. We give also some new results on the
integrability of Riemannian Lie algebroids.\\



\section{introduction}

Lie groupoids and Lie algebroids are now a central notion in
differential geometry and constitute an active domain of research.
They have many applications in various part of mathematics (see
for instance [4, 5, 6, 14]).
 Roughly speaking, a Lie algebroid is a structure where one
replaces the tangent bundle with a new vector bundle with similar
properties. In this spirit, many  geometrical notions which
involves  the tangent bundle were  generalized to the context of
Lie algebroids. For instance, covariant derivatives were
generalized by Fernandes [9], Lagrangian mechanics were generalized by Weinstein
[18] (see also [6]). Actually, a Riemannian metric on a manifold is a notion
which involves the Lie algebroid structure of the tangent bundle
and the Koszul formula, which defines the Levi-Civita connection,
is an illustration of this fact. A Riemannian metric on a Lie
algebroid is a natural extension of the classical notion of
Riemannian metric on a manifold and this notion appeared first
in the context of Lie algebroids associated to Poisson structures
(see [2, 3, 12, 13]).

In this paper, we present the basic concepts related to a
Riemannian structure on a Lie algebroid, namely, we will show that
most of the classical tools and results known in Riemannian
geometry can be stated in this setting after some slight
arrangements. In Section 2, we present some basic facts on
connections on  Lie algebroids based on recent results of [7]. In
Section 3, we define the Levi-Civita connection associated to a
Riemannian Lie algebroid and we show that this connection
satisfies relations which are similar to those introduced by
O'Neill in the context of Riemannian submersions [15] (see [1] for
a detailed presentation). Section 4 is devoted to the study of the
geodesic flow of a Riemannian Lie algebroid. As the classical
case, we define the analogous of the Sasaki metric and we compute
the divergence of the geodesic flow with respect to this metric. This
divergence does not vanish in general  contrast to  the classical
Liouville theorem. We state the first and the second variation
formulas and introduce the analogous of Jacobi fields. This
section can be thought of as a completion of subsection 4.2 in
[18] and Section 5 in [11]. In Section 5, we study the curvature
of a Riemannian Lie algebroid and generalize some classical
results, namely, Mayers  theorem. Section 6 is
devoted to the study of integrability of Riemannian Lie
algebroids, for instance, we show that a complete  Riemannian Lie algebroid
whose Riemannian curvature is nonpositive is integrable and is diffeomorphic to its Weinsein Lie groupoid. This is a generalization of Hadamard-Cartan theorem.

\section{Background  on Lie algebroids}

In this section we review some basic facts related to Lie
algebroids and to  connections in the context of  Lie algebroids
(see [6, 7, 9] for a detailed presentation).

\subsection{Canonical Poisson structure on the dual of a Lie algebroid}

A Lie algebroid $A$ over a smooth manifold $M$ is a vector bundle
$p:A\too M$ together with a Lie algebra structure $\br$ on the
space of sections $\Ga(A)$ and a bundle map $\#:A\too TM$ called
{\it anchor} such that

$(i)$ the induced map $\#:\Ga(A)\too\X(M)$ is a Lie algebra
homomorphism;

$(ii)$ for any sections $a,b\in\Ga(A)$ and for every smooth
function $f\in\C^\infty(M)$ we have the Leibniz identity
\begin{equation}[a,fb]=f[a,b]+\#(a)(f)b.\end{equation}
An immediate consequence of this definition is that, for any $x\in
M$, there is an induced Lie bracket say $\br_x$ on
$${\mathcal G}_x={\mathrm Ker}(\#_x)\subset A_x$$which makes it
into a Lie algebra.

The following theorem describes the local structure of a Lie
algebroid (for a proof see [9]). We denote by $n$ the dimension of
$M$ and by $r$ the rank of the vector bundle $A\too M$.

\begin{th} {\bf(Local splitting)} Let $x_0\in M$ be a point where
$\#_{x_0}$ has rank $q$. There exists a coordinates
$(x_1,\ldots,x_q,y_1,\ldots,y_{n-q})$ valid in a neighborhood $U$
of $x_0$ and a basis of sections $\{a_1,\ldots,a_r\}$ of $A$ over
$U$, such that
\begin{eqnarray*}
\#(a_i)&=&{\partial_{x_i}}\qquad(i=1,\ldots,q),\\
\#(a_i)&=&\sum_{j}b^{ij}{\partial_{y_j}}\qquad(i=q+1,\ldots,r),\end{eqnarray*}where
$b^{ij}\in C^\infty(U)$ are smooth functions depending only on the
$y'$s and vanishing at $x_0$: $b^{ij}=b^{ij}(y^s),$
$b^{ij}(x_0)=0$. Moreover, for any $i,j=1,\ldots,r$,
$$[a_i,a_j]=\sum_{u} C_{ij}^ua_u,$$where $C_{ij}^u\in C^\infty(U)$
vanish if $u\leq q$ and satisfy
$$\sum_{u>q}\frac{\partial C_{ij}^u}{\partial
x_s}b^{ut}=0.$$\end{th}

From this theorem we deduce that the image of $\#$ defines a
smooth generalized distribution in $M$, in the sense of Sussman
$[16]$, which is integrable. This foliation is called {\it
characteristic foliation} of $A$.  We call $A$  {\it transitive Lie algebroid } if $\#$ is
surjective, so the leaves are the connected components of $M$.
\\ We denote by $A_L$ the restriction of $A$ to a leaf $L$. From
$(1)$ one can deduce easily that the bracket $\br$ induces a
bracket on the space of sections of $p_L:A_L\too L$ and hence a
transitive Lie algebroid structure. When $x$ run over $L$ the
$\G_x'$s are all isomorphic and fit into a Lie algebra bundle
$\G_L$ over $L$ (see [14]).\\
The dual $A^*$ of a Lie algebroid  $p:A\too M$ carries a natural
Poisson structure which can be described as follows.

For any function $f\in C^\infty(A^*)$ and for any section $\xi\in\Ga(A^*)$, we define a section $f_\xi\in\Ga(A)$ by putting, for any $x\in M$ and for any $\mu_x\in A_x^*$,
$$<\mu_x,f_\xi(x)>=\frac{d}{dt}_{|t=0}f(\xi(x)+t\mu_x).$$Now, for any functions $f,g\in C^\infty(A^*)$, we define the bracket $\{f,g\}$ by putting, for any section $\xi\in\Ga(A^*)$,
\begin{equation}\{f,g\}\circ\xi=<\xi,[f_\xi,g_\xi]>+\#(f_\xi)(g\circ\sigma)-\#(g_\xi)(f\circ\sigma),\end{equation}
where $\sigma:M\too A^*$ is the zero section. One checks
that this bracket  defines a Poisson structure.

 If one chooses local coordinates $(x_1,\ldots,x_n)$
over a  neighborhood $U$ of $M$ and a basis of local
sections $(a_1,\ldots,a_r)$ over $U$, we have structure functions
$b^{si},C_{st}^u\in C^\infty(U)$ defined by
\begin{eqnarray*}
\#(a_s)&=&\sum_{i=1}^nb^{si}{\partial_{x_i}}\qquad(s=1,\ldots,r),\\
\;[a_s,a_t]&=&\sum_{u=1}^rC_{st}^ua_u\qquad(s,t=1,\ldots,r).\end{eqnarray*}
 Let
$(\xi_1,\ldots,\xi_r)$ denote the linear coordinates on the fibers
of $A^*$ associated with the dual basis $(a^1,\ldots,a^r)$.
One can see easily that
\begin{equation}\{x_i,x_j\}=0,\quad\{x_i,\xi_s\}=-b^{si}\quad\mbox{and}\quad
\{\xi_s,\xi_t\}=\sum_{u}C_{st}^u\xi_u.\end{equation}
\begin{exem}\begin{enumerate}\item The basic example of a Lie algebroid over $M$ is the tangent bundle itself, with the identity mapping as anchor. The associated Poisson structure on $T^*M$ is the one defined by the symplectic form $d\la$ where $\la$ is the Liouville form.
\item Every finite dimensional Lie algebra is a Lie algebroid over a one point space. The associated Poisson structure on the dual is the Lie-Poisson structure.
\item Any integrable subbundle of $TM$ is a Lie algebroid with the inclusion as anchor and the induced bracket.\item Let $(P,\pi)$ be a Poisson manifold. Then there is a natural Lie algebra structure on $\Om^1(P)$ which makes $T^*P$ into a Lie algebroid over $P$ (see [17]).

\end{enumerate}

\end{exem}

\subsection{Connections on Lie algebroids}
We develop now the basic theory of connections on Lie algebroids.
This notion, which is the natural extension of the usual concept
of covariant connection, have recently turned out to be useful in
the study of Lie algebroids. It  appeared first in the context of
Poisson geometry (see [9, 10, 17]).

Let $p:A\too M$ be a Lie algebroid with anchor map $\#$. An $A$-connection on a vector bundle $E\too M$  is
an operator $\na:\Ga(A)\times\Ga(E)\too\Ga(E)$ which satisfies:
\begin{enumerate}

 \item $\na_{a+b}s=\na_a s+\na_b s$ for any
$a,b\in\Ga(A)$ and $s\in\Ga(E)$;
 \item $\na_a(s_1+s_2)=\na_a s_1+\na_a s_2$ for any
$a\in\Ga(A)$ and $s_1,s_2\in\Ga(E)$;

 \item $\na_{fa}s=f\na_a s$ for any
$a\in\Ga(A)$, $s\in\Ga(E)$ and $f\in C^\infty(M)$;

 \item $\na_a(fs)=f\na_a s+\#(a)(f)s$
 for any
$a\in\Ga(A)$, $s\in\Ga(E)$ and $f\in C^\infty(M).$\end{enumerate}

From the definition, one can deduce immediately that, for any leaf
$L$, $\na$ induces an $A_L$-connection on $E_L\too L$.

 Given an
$A$-connection on a vector bundle $E$ over $M$, most of the
classical constructions (related to a classical covariant
derivative) extend to Lie algebroids, provided we use the
appropriate notion of paths on $A$.

\begin{Def} Let $p:A\too M$ be a Lie algebroid with anchor $\#$.\begin{enumerate}\item
 An
$A$-path is a smooth path $\al:[t_0,t_1]\too A$ such that $$
\#(\al(t))=\frac{d}{dt}p(\al(t)),\qquad\qquad t\in[t_0,t_1].$$ The
curve $\ga:[t_0,t_1]\too M$ given by $\ga(t)=p(\al(t))$ will be called
the {\it base path of } $\al$.\item An $A$-path $\al$ is called
vertical if $\al(t)\in\G_{p(\al(t_0))}$ for any
$t\in[t_0,t_1]$.\end{enumerate}\end{Def}

\begin{rem} Even if, for a vertical $A$-path the base path is
reduced to a constant curve, vertical $A$-paths play a non trivial
role in the study of connections on a Lie algebroid.\end{rem}

\subsection{Parallel transport} Let $p:A\too M$ be a Lie
algebroid, $E\too M$ a vector bundle  and $\na$  an $A$-connection
on $E$. Fix  an  $A$-path $\al:[t_0,t_1]\too A$. An $\al$-section of
$E$ is a smooth map  $s:[t_0,t_1]\too E$ such that the projections on
$M$ of $\al$ and $s$ define the same base path. We denote by
$\Ga(E)_\al$ the space of $\al$-sections of $E$. Then there is
exists an unique map
$$\na^\al:\Ga(E)_\al\too\Ga(E)_\al$$satisfying:

\begin{enumerate}

\item $\na^\al (c_1s_1+c_2s_2)=c_1\na^\al s_1+c_2\na^\al s_2$,
$c_1,c_2\in\reel$;

\item $\na^\al fs=f's+f\na^\al s$ where $f:[t_0,t_1]\too\reel$ is a
smooth function;

\item if $\wi s$ is a local section of $E$ which extends $s$  and
$\#(\al(t))\not=0$ then
$$\na^\al s(t)=\na_{\al(t)}\wi s;$$
\item if $\wi s$ is a local section of $E$ which extends $s$ and
$\al$ is vertical  then
$$\na^\al s(t)=\na_{\al(t)}\wi s+\frac{d}{dt}s(t).$$\end{enumerate}

An $\al$-section $s$ is called parallel along $\al$ if $\na^\al
s=0$.  One has then the notion of parallel transport along $\al$,
denoted by
$$\tau_\al^t:E_{\ga(t_0)}\too E_{\ga(t)},$$and
$\tau_\al^t(s_0)=s(t)$ where $s$ is the unique parallel
$\al$-section satisfying $s(0)=s_0$.

If $\al_0\in A_x$ and $s$ is a section of $E$ in a neighborhood of
$x$, one can check easily that
\begin{equation}\na_{\al_0}s=\frac{d}{dt}_{|t=0}(\tau_\al^t)^{-1}(s(\ga(t))),
\end{equation}where $\al$ is any $A$-path satisfying
$\al(0)=\al_0.$

\subsection{ Linear $A$-connections, geodesics and compatibility with the Lie algebroid structure} Let
$p:A\too M$ be a Lie algebroid with anchor $\#$.
We shall call  $A$-connections on the vector bundle $A\too M$
 {\it linear $A$-connections}.\\
Let $\D$ be a linear $A$-connection.  An $A$-path $\al:[t_0,t_1]\too
A$ is  a {\it geodesic} of $\D$ if $\D^\al\al=0$. Let $(x_1,\ldots,x_n)$
be a local system of coordinates on an open set $U$ and
$(a_1,\ldots,a_r)$ a basis of local sections over $U$. The structure functions $b^{si},C_{st}^u\in
C^\infty(U)$ are given by
\begin{eqnarray*}
\#a_s&=&\sum_{i=1}^nb^{si}\partial_{x_i}\qquad(s=1,\ldots,r),\\
\;[a_s,a_t]&=&\sum_{u=1}^rC_{st}^ua_u\qquad(s,t=1,\ldots,r).\end{eqnarray*}We
define the Christoffel symbols of $\D$ according to
$(a_1,\ldots,a_r)$ as usually by
$$\D_{a_s}a_t=\sum_{u=1}^r\Ga_{st}^ua_u.$$
 An $A$-path $\al$ is locally determined  by
$$\al(t)=\sum_{i=1}^r\al_i(t)a_i,\quad p(\al(t))=(x_1(t),\ldots,x_n(t)).$$
 The $A$-path $\al:[t_0,t_1]\too A$ is a geodesic if, for
$i=1,\ldots,n$ and $j=1,\ldots,r$,
\begin{equation}\left\{\begin{array}{ccc}
\dot
x_i(t)&=&\di\sum_{j=1}^r\al_j(t)b^{ji}(x_1(t),\ldots,x_n(t)),\\
\dot
\al_j(t)&=&-\di\sum_{s,u=1}^r\al_s(t)\al_u(t)\Ga_{su}^j(x_1(t),\ldots,x_n(t)).\end{array}\right.\end{equation}

Exactly as in the classical case, one has existence and uniqueness
of geodesics with given initial base point $x\in M$ and "initial
speed" $a_0\in A_x$. Actually, there exists a vector field $G$ on $A$ such that the geodesics of $\D$ are the integral curves of $G$. We call $G$ the {\it geodesic vector field} associated to $\D$ and $\D$ is called complete if $G$ is complete.\bigskip

We introduce now two natural notions of compatibility between a
linear $A$-connection and the structure of Lie algebroid.
\begin{Def}
\begin{enumerate}\item A linear $A$-connection $\D$ is strongly
compatible with the Lie algebroid structure  if, for any $A$-path
$\al$ , the parallel transport  $\tau_\al$ preserves ${\mathrm
Ker}\#$. \item A linear $A$-connection $\D$ is weakly compatible
with the Lie algebroid structure  if, for any vertical $A$-path
$\al$, the parallel transport $\tau_\al$ preserves ${\mathrm
Ker}\#$. \end{enumerate}\end{Def}

The following proposition gives an useful characterization of the
two notions of compatibility above.

\begin{pr}
\begin{enumerate}\item A linear $A$-connection $\D$ is strongly
compatible with the Lie algebroid structure  if and only if, for
any leaf $L$ any sections $\al\in\Ga(A_L)$ and $\be\in\Ga(\G_L)$,
$\D_\al\be\in\Ga(\G_L)$. \item A linear $A$-connection $\D$ is
weakly compatible with the Lie algebroid structure  if and only
if, for any leaf $L$ and for any sections $\al\in\Ga(\G_L)$ and
$\be\in\Ga(\G_L)$, $\D_\al\be\in\Ga(\G_L)$.

\end{enumerate}\end{pr}

{\bf Proof.} This is a consequence of $(4)$.$\Box$

\begin{exem}Let $p:A\too M$ be a Lie algebroid and  $\na$ be a  $TM$-connection on $A$. Associated with $\na$ there
is an obvious linear $A$-connection
$$\D^0_ab=\na_{\#(a)}b$$ which is clearly weakly compatible with
the Lie algebroid structure. A bit more subtle is the following
linear $A$-connection
$$\D^1_ab=\na_{\#(b)}a+[a,b]$$which is strongly compatible with
the Lie algebroid structure. These connections play a fundamental
role in the theory of characteristic classes (see for instance
[9]).\end{exem}

\begin{rem} In [9] there is a notion of compatibility between
linear $A$-connections and the Lie algebroid structure which is
stronger than the notion of compatibility given in Definition 1.2
1. \end{rem}

 \subsection{ Variations of $A$-paths, homotopy and curvature of $A$-connections}
We give an interpretation of the torsion and the curvature of an
$A$-connection which leads naturally to the notion of homotopy of
$A$-paths. This notion plays a crucial role in the integrability
of Lie algebroids (see [7]).

 Let $p:A\too M$ be a Lie
algebroid with anchor $\#$ and $E\too M$ a vector bundle. The
curvature of an $A$-connection $\na$ on $E$ is formally identical
to the usual definition
$$R(a,b)s=\na_{a}\na_{b}s-\na_{b}\na_{a}s-\na_{[a,b]}s,$$
where $a,b\in\Ga(A)$ and $s\in\Ga(E)$. The connection $\na$ is
called flat if $R$ vanishes identically.

If $\D$ is  a linear $A$-connection the torsion of $\D$ is given
by
$$T_\D(a,b)=\D_ab-\D_ba-[a,b].$$

In the usual case,  the curvature and the torsion can be
interpreted by using variations of paths. Let us precise this
well-known fact in our context. First, let us give the appropriate
notion of variation of paths.

A {\it variation} of $A$-paths is a smooth map
$\al:[0,1]\times[0,1]\too A$, $(\e,t)\mapsto \al(\e,t)$ such that:

$(i)$ for any $\e\in[0,1]$, the map $\al(\e,.)$ is an $A$-path,

$(ii)$ the base variation $\ga(\e,t)=p(\al(\e,t))$ lies entirely
in a fixed leaf $L$ of the characteristic foliation.

A variation of $A$-path $\al$ is given, we call a smooth map
$\be:[0,1]\times[0,1]\too A$ {\it transverse variation} to $\al$
if $\al$ and $\be$ have the same base variation $\ga$ and
$\#(\be)=\frac{\partial\ga}{\partial\e}$. \\ Fix $(\al,\be)$ as above and denote
by $\ga$ the commune base path. Let $\na$ be an  $A$-connection on
a vector bundle $E\too M$ and let $s:[0,1]\times [0,1]\too E$ be a
section over $\ga$. For any $\e\in [0,1]$, $t\mapsto\al(\e,t)$ is an
$A$-path and we denote by $\na_ts$ the derivative of $t\mapsto s(\e,t)$
along this $A$-path. On the other hand, for any $t\in [0,1]$,
$\e\mapsto\be(\e,t)$ is an $A$-path and we denote by $\na_\e s$ the
derivative of $\e\mapsto s(\e,t)$ along this $A$-path.\\
It is clear that if $\#$ is injective, there is an unique
transverse variation to a given variation of $A$-paths. However,
if $\#$ is not injective, a given variation of $A$-paths  admits
many transverse variations to it. There is a way which permit the
control  of transverse variations
 to a fixed variation of $A$-path. Let us explain this important
 fact which is at the origin of the  notion of homotopy
 of $A$-paths used in [7].
The first claim in the following proposition is a reformulation of
a part of Proposition 1.3 in [7].

\begin{pr} With the notation above the following assertions hold.
\begin{enumerate}

\item For any linear $A$-connection $\D$, the variation
$$\De(\al,\be)=\D_t\be-\D_\e \al- T_\D(\al,\be)$$ does
not depend on $\D$ and satisfies $\#(\De(\al,\be))=0$. \item for
any $A$-connection $\na$ on $E$ and for any section $s$ of $E$
over $\ga$ $$ \na_t\na_\e
s-\na_\e\na_ts=R(\al,\be)s+\na_{\De(\al,\be)}s.$$
\end{enumerate}\end{pr}

{\bf Proof.}\begin{enumerate}\item Fix $(\e_0,t_0)\in [0,1]\times
[0,1]$ and choose a local coordinates
$(x_1,\ldots,x_q,y_1,\ldots,y_{n-q})$  near $x_0=\ga(\e_0,t_0)$
and a basis of sections $(a_1,\ldots,a_r)$ as in Theorem 1.1
($q={\mathrm rank}\#_{x_0}$). In these coordinates, we have
\begin{equation}
\left\{\begin{array}{l}\di \al(\e,t)=\sum_{i=1}^r\al^i(\e,t)a_i,\\
\di \be(\e,t)=\sum_{i=1}^r\be^i(\e,t)a_i,\\
\ga(\e,t)=(x_1(\e,t),\ldots,x_q(\e,t),c_1,\ldots,c_{n-q}),\\
\di \frac{\partial\ga}{\partial t}=\sum_{j=1}^q\frac{\partial x_j}{\partial
t}{\partial}_{x_j}=
\sum_{i=1}^q\al^j(\e,t){\partial}_{x_j},\\
\di \frac{\partial\ga}{\partial\e}=\sum_{j=1}^q\frac{\partial
x_j}{\partial\e}{\partial}_{
x_j}=\sum_{i=1}^q\be^j(\e,t){\partial}_{ x_j},
\end{array}\right.\end{equation}where $c_1,\ldots,c_{n-q}$ are
constant. Now $$\D_t \be=\sum_{i=1}^r\frac{\partial
\be^i}{\partial
t}a_i+\sum_{i,j=1}^r\al^j\be^i\D_{a_j}a_i\quad\mbox{and}\quad
\D_\e
\al=\sum_{i=1}^r\frac{\partial\al^i}{\partial\e}a_i+\sum_{i,j=1}^r\al^i\be^j
\D_{a_j}a_i.$$

Hence $$\D_t \be-\D_\e \al=\sum_{i=1}^r\left(\frac{\partial
\be^i}{\partial t}-\frac{\partial \al^i}{\partial
\e}\right)a_i+T_\D(\al,\be)+\sum_{i,j=1}^r\al^i\be^j[{a_i},a_j].$$
Now, form $(6)$, we have $\di\frac{\partial \be^i}{\partial
t}=\frac{\partial \al^i}{\partial \e}$ for any $i=1,\ldots,q$, so
\begin{equation}\D_t \be-\D_\e
\al-T_\D(\al,\be)=\sum_{i=q+1}^r\left(\frac{\partial
\be^i}{\partial t}-\frac{\partial \al^i}{\partial
\e}\right)a_i+\sum_{i,j=1}^r\al^i\be^j[{a_i},a_j].\end{equation}

One can see that the right hand of this equality lies in ${\mathrm
Ker}\#$ and does not depend on $\D$.

\item We choose a local trivialization
$(x_1,\ldots,x_q,y_1,\ldots,y_{n-q},a_1,\ldots,a_r)$ as above, we
trivialize $E$ near $x_0$ by a local basis of sections
$(e_1,\ldots,e_\mu)$ and put
$$s(\e,t)=\sum_{j=1}^\mu s^j(\e,t)e_j.$$

We have
\begin{eqnarray*}
\na_ts&=&\sum_{j=1}^\mu\frac{\partial s^j}{\partial
t}e_j+\sum_{i,j}\al^is^j\na_{a_i}e_j.\\
\na_\e\na_ts&=&\sum_{j=1}^\mu\frac{\partial^2 s^j}{\partial
\e\partial t}e_j+\sum_{i,j}\left(\be^i\frac{\partial s^j}{\partial
t}+\frac{\partial \al^i}{\partial \e}s^j+\al^i\frac{\partial
s^j}{\partial
\e}\right)\na_{a_i}e_j+\sum_{i,j,k}\be^k\al^is^j\na_{a_k}\na_{a_i}e_j.\\
\na_t\na_\e s&=&\sum_{j=1}^\mu\frac{\partial^2 s^j}{\partial
t\partial \e}e_j+\sum_{i,j}\left(\al^i\frac{\partial s^j}{\partial
\e}+\frac{\partial \be^i}{\partial t}s^j+\be^i\frac{\partial
s^j}{\partial
t}\right)\na_{a_i}e_j+\sum_{i,j,k}\al^k\be^is^j\na_{a_k}\na_{a_i}e_j.\\
\na_t\na_\e s&-&\na_\e\na_t
s-R(\al,\be)s=\sum_{i,j}\left(\frac{\partial \be^i}{\partial
t}-\frac{\partial \al^i}{\partial
\e}\right)s^j\na_{a_i}e_j+\sum_{i,j,k}\al^k\be^is^j\na_{[a_k,a_i]}e_j.\end{eqnarray*}
The above computation and $(7)$ give the desired
formula.$\Box$\end{enumerate}

From the expression of $\De(\al,\be)$ given by $(7)$ and from
$(6)$, we have
\begin{equation}\De(\al,\be)=0\Leftrightarrow\left\{\begin{array}{lc}\di
\frac{\partial\al_i}{\partial \e}-\frac{\partial\be_i}{\partial
t}=\sum_{l,k=1}^r\al^l\be^kC_{lk}^i&\di i=q+1,\ldots,r,\\\di
\di\al^j=\frac{\partial x_j}{\partial t},\;\be^j=\frac{\partial
x_j}{\partial \e}& j=1,\ldots,q.\end{array}\right.\end{equation}

 Now by using  the standard results about linear
differential systems one can deduce easily the following useful
proposition (compare to Proposition 1.1 in [7]).

\begin{pr} Let $p:A\too M$ be a Lie algebroid. Then, for a given
variation of $A$-paths $\al$ and for given $\be_0:[0,1]\too A$
such that $\#(\be_0)(\e)=\frac{\partial p\circ\al}{\partial\e}(\e,0)$ there exists an
unique transverse variation $\be$ to $\al$ such that
$$\De(\al,\be)=0\quad\mbox{and}\quad \be(\e,0)=\be_0(\e)\quad\mbox{for any}\;\e\in[0,1].$$
\end{pr}

Following [7], we can now define the homotpoy  of $A$-paths with
fixed end-points. Let $\al_0$ and $\al_1$ be two $A$-paths on a
Lie algebroid $p:A\too M$ such that $p(\al_0(0))=p(\al_1(0))$ and
$p(\al_0(1))=p(\al_1(1))$. An $A$-homotopy with fixed end-points from
$\al_0$  to $\al_1$ is a variation of $A$-paths $\al$ such that:

$(i)$ $p(\al(\e,0))=p(\al(0,0))$  and $p(\al(\e,1))=p(\al(0,1))$
for any $\e\in[0,1]$, $\al(0,.)=\al_0$ and $\al(1,.)=\al_1$,

$(ii)$ the unique transverse variation $\be$ to $\al$ satisfying
$\De(\al,\be)=0$ and $\be(\e,0)=0$ satisfies also $\be(\e,1)=0$.

The following Lemma will be useful latter.

\begin{Le} Let $\al_0:[0,1]\too A$ be an  $A$-path and
$\be_0:[0,1]\too A$ an $\al_0$-section  such that
$\be_0(0)=\be_0(1)=0$. Then there exists an $A$-homotopy $\al$ with fixed
end-points  such that $\al(0,.)=\al_0$ and the
corresponding transverse variation $\be$ satisfies
$\be(0,.)=\be_0$.\end{Le}

{\bf Proof.} Consider the base path $\ga_0:[0,1]\too M$ of $\al_0$
and choose an homotopy $\ga:[0,1]\times [0,1]\too M$ with fixed
end points such that $\ga$ lies in the same leaf as $\ga_0$,
$\ga(0,.)=\ga_0$ and
$\frac{\partial\ga}{\partial\e}(0,t)=\#(\be_0(t))$. We choose also
$\be:[0,1]\times [0,1]\too A$ such that $\be(0,t)=\be_0(t)$ for
any $t\in[0,1]$, $\be(\e,0)=\be(\e,1)=0$ for any $\e\in[0,1]$ and
$\frac{\partial\ga}{\partial\e}(\e,t)=\#(\be(\e,t))$ for any
$(\e,t)$. From $(8)$, one can deduce that there exists an unique
variation $\al:[0,1]\times[0,1]\too A$ such that the base path of
$\al$ is $\ga$, $\frac{\partial\ga}{\partial
t}(\e,t)=\#(\al(\e,t))$, $\al(0,.)=\al_0$ and $\De(\al,\be)=0.$
This variation  is clearly an $A$-homotopy with fixed end-points and
satisfies the required properties.$\Box$

\section{Riemannian metrics on Lie algebroids}

In this section, we  introduce the notion of Riemannian metric on
a Lie algebroid which is a natural extension of the notion of
Riemannian metric on a manifold. We  show that most of the
classical notions associated to a Riemannian metric can be defined
in this context, namely, Levi-Civita connection, geodesics,
geodesic flow, Sasaki metric, first and second variation formulas,
Jacobi fields, the exponential... We  show also that the
Riemannian curvature of a Riemannian metric on a Lie algebroid
satisfies formulas which are formally identical to the O'Neill
formulas for Riemannian submersions.

\subsection{The Levi-Civita connection of a Riemannian metric on a
Lie algebroid}

A Riemannian metric on a Lie algebroid $p:A\too M$ is the data,
for any $x\in M$, of  a scalar product $\prs_x$ on the fiber $A_x$
such that, for any local section $a,b\in\Ga(A)$, the function
$<a,b>$ is smooth.

The most interesting fact about Riemannian metrics on Lie
algebroids is the existence on the analogous of the Levi-Civita
connection. Indeed, if $\prs$ is a Riemannian metric on a Lie
algebroid $p:A\too M$, then the formula
\begin{eqnarray*}
2<\D_ab,c>&=&{\#}(a).<b,c>+{\#}(b).<a,c>-
{\#}(c).<a,b>\\
&+&<[c,a],b>+<[c,b],a>+<[a,b],c>\end{eqnarray*}defines
 a linear $A$-connection  which is characterized by the two following properties:

$(i)$ $\D$ is metric, i.e., $\#(a).<b,c>=<\D_ab,c>+<b,\D_ac>$,

$(ii)$ $\D$ is torsion free, i.e., $\D_ab-\D_ba=[a,b].$

We  call $\D$ the {\it Levi-Civita $A$-connection} associated to
the Riemannian metric $\prs$. \\
In  local coordinates $(x_1,\ldots,x_n)$ over a trivializing
neighborhood $U$ of $M$ where $A$ admits a basis of local sections
$(a_1,\ldots,a_r)$, the Levi-Civita $A$-connection is determined
by the Christoffel's symbols defined by
$\D_{a_i}a_j=\sum_{k}\Ga_{ij}^ka_k$. We have \begin{eqnarray}
\Ga_{ij}^k&=&\frac12\sum_{l=1}^r\sum_{u=1}^ng^{kl}\left(b^{iu}\partial_{x_u}(g_{jl})+
b^{ju}\partial_{x_u}(g_{il})-b^{lu}\partial_{x_u}(g_{ij})\right)\nonumber\\
&&+\frac12\sum_{l=1}^r\sum_{u=1}^rg^{kl}\left(
C_{ij}^ug_{ul}+C_{li}^ug_{uj}+C_{lj}^ug_{ui}\right)\end{eqnarray}
where the structure functions $b^{si},C_{st}^u\in C^\infty(U)$ are
given by
\begin{eqnarray*}
\#a_s&=&\sum_{i=1}^nb^{si}{\partial_{x_i}}\qquad(s=1,\ldots,r),\\
\;[a_s,a_t]&=&\sum_{u=1}^rC_{st}^ua_u\qquad(s,t=1,\ldots,r),\end{eqnarray*}
 $g_{ij}=<a_i,a_j>$ and  $(g^{ij})$ denotes the inverse
matrix of $(g_{ij})$.\\
As the classical case, for any
$A$-path $\al$ and for any $\al$-sections $\be$ and $\ga$, one
has
\begin{equation}\frac{d}{dt}<\be,\ga>=<\D^\al\be,\ga>+<\be,\D^\al\ga>.\end{equation}
\begin{rem}There are two extremal cases:
\begin{enumerate}\item The Lie algebroid $A$ is the tangent bundle $TM$ of a
manifold and we recover the classical notion of Riemannian
manifold. \item The Lie algebroid $A$ is a Lie algebra $\G$
considered as a Lie algebroid over a point. In this case a
Riemannian metric on $\G$ is a scalar product $\prs$ and the
Levi-Civita $\G$-connection is the product $\D:\G\times\G\too\G$
given by
$$2<\D_uv,w>=<[u,v],w>+<[w,u],v>+<[w,v],u>.$$ Actually $\D$ is the
infinitesimal data associated to the Levi-Civita connection of the
left invariant metric associated to $\prs$ on any Lie group with
$\G$ as a Lie algebra.\end{enumerate}\end{rem}

The general setting is a combination of these two extremal cases. Indeed,
let $\prs$ be a Riemannian metric on a Lie algebroid $p:A\too M$
with anchor $\#$, then we have:
\begin{enumerate}\item
For any leaf $L$  of the
characteristic foliation and  for any $x\in L$, we have
$$A_x=\G_x\oplus\G_x^{\perp},$$where $\G_x^{\perp}$ is the
orthogonal of $\G_x$ with respect $\prs_x$. The restriction of the
anchor $\#$ to $\G_x^{\perp}$ is an isomorphism into $T_xL$ and
hence induces  a scalar product on $T_xL$
$$<u,v>_L=<a,b>,$$ where $a,b\in\G_x^{\perp}$ and $\#(a)=u$ and
$\#(v)=b$. Thus $\prs$ induces a Riemannian metric $\prs_L$ on
$L$. We call it the {\it induced Riemannian metric} on $L$.\item
The scalar product $\prs_x$ induce a scalar
product on $\G_x$ and we denote by $\widehat{\D}$ the Levi-Civita
$\G_x$-connection associated with $(\G_x,\prs_x)$.\end{enumerate}

Let us  precise more this situation. Fix a leaf $L$ and consider $p_L:A_L\too L$. We have
$$A_L=\G_L\oplus\G_L^\perp.$$
We call the elements of $\Ga(\G_L)$ {\it vertical sections} and the
elements of $\Ga(\G_L^\perp)$ {\it horizontal sections}. For any section
$a$, we denote by $a^v$ its vertical component and by $a^h$ its
horizontal component. Note that the bracket of a vertical section with every section is a vertical section. Thus, in the Riemannian point of view, the short exact sequence
$$0\too\G_L\too A_L\too TL$$ is formally identical to a Riemannian
submersion. So we can introduce the O'Neill tensors [15] (see [1]
for a detailed presentation).

We denote by $T$ and $H$ the elements of $\Ga(A\otimes A\otimes A^*)$ whose values on sections
$a,b$ are given by
$$T_ab=(\D_{a^v}b^v)^h+(\D_{a^v}b^h)^v\quad\mbox{and}\quad H_ab=(\D_{a^h}b^v)^h+(\D_{a^h}b^h)^v.$$

The following properties of $T$ and $H$ are easy consequence of the
definition.
\begin{eqnarray*}
T_{a^h}b^v&=&T_{a^h}b^h=0,\\
T_{a^v}b^v&=&(\D_{a^v}b^v)^h\quad\mbox{and}\quad
T_{a^v}b^h=(\D_{a^v}b^h)^v,\\
T_{a^v}b^v&=&T_{b^v}a^v,\\
<T_{a^v}b^v,c^h>&=&-<T_{a^v}c^h,b^v>,
\\H_{a^v}b^h&=&H_{a^v}b^v=0,\\
H_{a^h}b^v&=&(\D_{a^h}b^v)^h\quad\mbox{and}\quad
H_{a^h}b^h=(\D_{a^h}b^h)^v,\\
H_{a^h}b^h&=&-H_{b^h}a^h,\\
<H_{a^h}b^h,c^v>&=&-<H_{a^h}c^v,b^h>.\end{eqnarray*}

On the other hand
\begin{equation}H_{a^h}b^h=\frac12[a^h,b^h]^v,\end{equation}
and, for any $u,v\in\G_x$,
\begin{equation} \D_uv=\widehat{\D}_uv+T_uv.\end{equation}
Moreover, we have, for any $a,b\in\Ga(A)$,
\begin{eqnarray*}
\D_{a^v}b^h&=&T_{a^v}b^h+(\D_{a^v}b^h)^h,\\
\D_{a^h}b^v&=&(\D_{a^h}b^v)^v+H_{a^h}b^v,\\
\D_{a^h}b^h&=&H_{a^h}b^h+(\D_{a^h}b^h)^h.\end{eqnarray*}
The following proposition is an immediate consequence of the last relation.
\begin{pr} Let $\ga:[t_0,t_1]\too L$ be a smooth path and let $\ga^h:[t_0,t_1]\too\G_L^\perp$ be the unique $A$-path with the base path  $\ga$. Then $\ga$ is a geodesic with respect to the induced Riemannian metric on $L$ if and only if $\ga^h$ is a geodesic of the Levi-Civita $A$-connexion.\end{pr}

The following proposition gives an interpretation of the tensors
$T$ and $H$.

\begin{pr} \begin{enumerate}\item The Levi-Civita $A$-connection
is strongly compatible with the Lie algebroid structure if and
only if $T=H=0.$\item The Levi-Civita $A$-connection is weakly
compatible with the Lie algebroid structure if and only if $T=0.$
\end{enumerate}\end{pr}

{\bf Proof.}  This is a consequence of
Proposition 1.1 and the relations above.$\Box$

\subsection{Geodesic flow of a Riemannian  Lie algebroid}

The Riemannian structure on a Lie algebroid $A$ gives arise to an identification between $A$ and its dual $A^*$. Thus $A$ inherits  a Poisson structure from the canonical Poisson structure of $A^*$. As the classical case (when $A=TM$), the Hamiltonian vector field associated to the energy function on $A$ is the geodesic flow of the Riemannian Lie algebroid. In this section, we give a complete proof of this fact and we generalize all the classical notions related to the geodesic flow, namely, the Sasaki metric, the first and second variation formulas, the Jacobi fields and the exponential. We give also the explicit formula of the divergence of the geodesic flow according to the Sasaki metric and we point out natural obstructions to the vanishing of this divergence. These obstructions vanish when $A$ is the tangent bundle of a manifold and we recover the classical Liouville theorem.\\

Let $p:A\too M$ be a Lie algebroid and $\prs$ a Riemannian metric
on $A$. The Riemannian metric defines a bundle  isomorphism
between $A$ and $A^*$ which transport the Lie-Poisson structure on
$A^*$ into a Poisson structure say $\pi_{\prs}$ in $A$.  Let
$E:A\too\reel$ be the energy function given
by$$E(a)=\frac12<a,a>$$ and let $X_E$ denote the hamiltonian vector field associated to $E$ with respect to
$\pi_{\prs}$.
The following result is a generalization
of a well-known result in Riemannian geometry.
\begin{th} The geodesics of the Levi-Civita connection associated to
$\prs$  are the integral curves of the hamiltonian vector field
$X_E$.\end{th}

{\bf Proof.} Let $(x_1,\ldots,x_n)$ be a system of coordinates over an open set $U$ of $M$ where $A$ admits a basis of
local sections $(a_1,\ldots,a_r)$. The structure functions
$b^{si},C_{st}^u\in C^\infty(U)$ are given by
\begin{eqnarray*}
\#a_s&=&\sum_{i=1}^nb^{si}{\partial_{x_i}}\qquad(s=1,\ldots,r),\\
\;[a_s,a_t]&=&\sum_{u=1}^rC_{st}^ua_u\qquad(s,t=1,\ldots,r).\end{eqnarray*}

We denote by $(\mu_1,\ldots,\mu_r)$ the linear coordinates on the
fibers of $A$ associated to $(a_1,\ldots,a_r)$ and by
$(\xi_1,\ldots,\xi_r)$ its dual coordinates on $A^*$. Recall that
the Poisson brackets on $A^*$ are given by
$$\{x_i,x_j\}=0,\quad\{x_i,\xi_s\}=-b^{si}\quad\mbox{and}\quad
\{\xi_s,\xi_t\}=\sum_{u}C_{st}^u\xi_u.$$

Put $g_{ij}=<a_i,a_j>$ and denote by $(g^{ij})$ the inverse
matrix of $(g_{ij})$. The isomorphism $\prs^\#:A^*\too A$, the
energy function and $X_E$ are  given, respectively , by
 \begin{eqnarray*}
 (x_1,\ldots,x_n,\xi^1,\ldots,\xi^r)&\mapsto&(x_1,\ldots,x_n,\sum_{i=1}^rg^{1i}\xi_i,\ldots,
\sum_{i=1}^rg^{ri}\xi_i),\\
E&=&\frac12\sum_{i,j}g_{ij}\mu_i\mu_j,\\
X_E&=&\sum_{i=1}^n\{E,x_i\}\partial_{x_i}+\sum_{j=1}^r\{E,\mu_j\}\partial_{\mu_j}.
\end{eqnarray*}

According to $(5)$, we must show that, for $i=1,\ldots,n$ and $j=1, \ldots,r$,
\begin{equation}
\{E,x_i\}=\sum_k\mu_kb^{ki}\quad\mbox{and}\quad
\{E,\mu_j\}=-\sum_{s,t}\mu_s\mu_t\Ga_{st}^j,\end{equation}where
$\Ga_{st}^j$ are the Christoffel symbols given by $(9)$, i.e.,\begin{eqnarray*}
\Ga_{ij}^k&=&\frac12\sum_{l=1}^r\sum_{u=1}^ng^{kl}\left(b^{iu}\partial_{x_u}(g_{jl})+
b^{ju}\partial_{x_u}(g_{il})-b^{lu}\partial_{x_u}(g_{ij})\right)\nonumber\\
&&+\frac12\sum_{l=1}^r\sum_{u=1}^rg^{kl}\left(
C_{ij}^ug_{ul}+C_{li}^ug_{uj}+C_{lj}^ug_{ui}\right).\end{eqnarray*}

\begin{enumerate}\item The first relation in $(13)$ is a straightforward computation. Indeed,
\begin{eqnarray*}
\{E,x_i\}&=&\frac12\sum_{k,l}g_{kl}\{\mu_k\mu_l,x_i\}
=\frac12\sum_{k,l}g_{kl}\left(\mu_k\{\mu_l,x_i\}+\mu_l\{\mu_k,x_i\}\right)\\
&=&\sum_{k,l}g_{kl}\mu_k\{\mu_l,x_i\}
=\sum_{k,l}g_{kl}\mu_k\{\sum_jg^{lj}\xi_j,x_i\}\\
&=&\sum_{k,l,j}g_{kl}g^{lj}\mu_k\{\xi_j,x_i\}
=\sum_{k,l,j}g_{kl}g^{lj}\mu_kb^{ji}\\
&=&\sum_{k,j}\left(\sum_lg_{kl}g^{lj}\right)\mu_kb^{ji}
=\sum_k\mu_kb^{ki}.\end{eqnarray*}

\item We must work much more  to establish the second relation in $(13)$.

Note first that
\begin{eqnarray}
 \sum_{s,t}\mu_s\mu_t\Ga_{st}^j&=&\frac12\sum_{s,t,u,l}g^{jl}\left(b^{su}\partial_{x_u}(g_{tl})+
b^{tu}\partial_{x_u}(g_{sl})-b^{lu}\partial_{x_u}(g_{st})\right)\mu_s\mu_t\nonumber\\
&&+\frac12\sum_{s,t,u,l}g^{jl}\left(
C_{st}^ug_{ul}+C_{ls}^ug_{ut}+C_{lt}^ug_{us}\right)\mu_s\mu_t\\
&\stackrel{(a)}=&\sum_{s,t,u,l}g^{jl}\left(b^{su}\partial_{x_u}(g_{tl})
-\frac12b^{lu}\partial_{x_u}(g_{st})
\right)\mu_s\mu_t+\sum_{s,t,u,l}g^{jl}g_{ut}
C_{ls}^u\mu_s\mu_t.\nonumber\end{eqnarray} We have used in $(a)$ the fact that $C_{st}^u=-C_{ts}^u$.  Now
\begin{eqnarray}
2\{E,\mu_j\}&=&\sum_{s,t}\{g_{st}\mu_s\mu_t,\mu_j\}\nonumber\\
&=&\sum_{s,t}\left(g_{st}\{\mu_s\mu_t,\mu_j\}+\mu_s\mu_t\{g_{st},\mu_j\}\right)\nonumber\\
&=&\sum_{s,t}\left(g_{st}\mu_s\{\mu_t,\mu_j\}+g_{st}\mu_t\{\mu_s,\mu_j\}\right)+
\sum_{s,t,l}\mu_s\mu_tg^{jl}\{g_{st},\xi_l\}\nonumber\\
&=&2\sum_{s,t}g_{st}\mu_s\{\mu_t,\mu_j\}-\sum_{s,t,l,u}
g^{jl}b^{lu}\partial_{x_u}(g_{st})\mu_s\mu_t.\end{eqnarray} By comparing $(14)$ and $(15)$, one can see that the desired relation is equivalent to
\begin{eqnarray}\sum_{s,t}g_{st}\mu_s\{\mu_t,\mu_j\}&=&
\sum_{s,t,u,l}g^{jl}\left(-b^{su}\partial_{x_u}(g_{tl})
+b^{lu}\partial_{x_u}(g_{st})
\right)\mu_s\mu_t\nonumber\\&&-\sum_{s,t,u,l}g^{jl}g_{ut}
C_{ls}^u\mu_s\mu_t.\end{eqnarray}Let us establish this relation. Note first that
\begin{eqnarray}
\{\mu_i,\mu_j\}&=&\sum_{k,l}\{g^{il}\xi_l,g^{jk}\xi_k\}\nonumber\\
&=&\sum_{k,l}\left(g^{il}g^{jk}\{\xi_l,\xi_k\}+g^{il}\xi_k\{\xi_l,g^{jk}\}+
g^{jk}\xi_l\{g^{il},\xi_k\}\right)\nonumber\\
&=&\sum_{k,l,u}g^{il}g^{jk}C_{lk}^u\xi_u+
\sum_{k,l,u}g^{il}\xi_kb^{lu}\partial_{x_u}(g^{jk})
-\sum_{k,l,u}g^{jk}\xi_lb^{ku}\partial_{x_u}(g^{il})\nonumber\\
&=&\sum_{k,l,u}g^{il}g^{jk}C_{lk}^u\xi_u+\sum_{k,l,u}b^{lk}\left(g^{il}\partial_{x_k}(g^{ju})
-g^{jl}\partial_{x_k}(g^{iu})\right)\xi_u.\nonumber
\end{eqnarray}
Hence
\begin{eqnarray*}
\sum_{s,t}g_{st}\mu_s\{\mu_t,\mu_j\}&=&\sum_{s,t,k,l,u}g_{st}g^{tl}g^{jk}C_{lk}^u\mu_s\xi_u\\
&&+\sum_{s,t,k,l,u}b^{lk}\left(g^{tl}\partial_{x_k}(g^{ju})
-g^{jl}\partial_{x_k}(g^{tu})\right)g_{st}\mu_s\xi_u.\end{eqnarray*}
Now
\begin{eqnarray*}
\sum_{s,t,k,l,u}g_{st}g^{tl}g^{jk}C_{lk}^u\mu_s\xi_u&=&
\sum_{s,k,u}g^{jk}C_{sk}^u\mu_s\xi_u
=\sum_{s,t,k,u}g^{jk}C_{sk}^u\mu_sg_{ut}\mu_t\\
&=&-\sum_{s,t,u,l}g^{jl}g_{ut}C_{ls}^u\mu_s\mu_t.\\
\sum_{s,t,k,l,u}b^{lk}g^{tl}\partial_{x_k}(g^{ju})
g_{st}\mu_s\xi_u&=&\sum_{s,k,u}b^{sk}\partial_{x_k}(g^{ju})
\mu_s\xi_u\\
&=&\sum_{s,t,k,u}b^{sk}g_{ut}\partial_{x_k}(g^{ju})
\mu_s\mu_t\\
&=&\sum_{s,t,k,u}b^{sk}\partial_{x_k}(g_{ut}g^{ju})
\mu_s\mu_t-\sum_{s,t,k,u}b^{sk}g^{ju}\partial_{x_k}(g_{ut})
\mu_s\mu_t\\
&=&-\sum_{s,t,u,l}b^{su}g^{jl}\partial_{x_u}(g_{lt})
\mu_s\mu_t.\\
\sum_{s,t,k,l,u}b^{lk}g^{jl}\partial_{x_k}(g^{tu})g_{st}\mu_s\xi_u&=&
-\sum_{s,t,k,l,u}b^{lk}g^{jl}\partial_{x_k}(g_{st})g^{tu}\mu_s\xi_u\\
&=&-\sum_{s,t,k,l,u,h}b^{lk}g^{jl}\partial_{x_k}(g_{st})g^{tu}g_{uh}\mu_s\mu_h\\
&=&-\sum_{s,t,k,l}b^{lk}g^{jl}\partial_{x_k}(g_{st})\mu_s\mu_t\\
&=&-\sum_{s,t,u,l}b^{lu}g^{jl}\partial_{x_u}(g_{st})\mu_s\mu_t.\end{eqnarray*}

Thus we get $(16)$ and the theorem follows.$\Box$\end{enumerate}

The flow of the Hamiltonian  vector field $X_E$ is called {\it the
geodesic flow of $\prs$}.
\begin{rem} Let $p:A\too M$ be a Riemannian Lie algebroid. Then:
\begin{enumerate}\item For any leaf $L$, the geodesic vector field $X_E$ is tangent to $A_L$ and to $\G_x$ for any $x\in L$. This follows from the fact that geodesics are $A$-paths.\item From Proposition 3.1, one can deduce that, for any leaf $L$, the geodesic vector field $X_E$ is tangent to $\G_L^\perp$.\end{enumerate}\end{rem}

\begin{co} Let $p:A\too M$ be Riemannian Lie algebroid. Then
\begin{enumerate}\item If $L$ is a compact leaf then the geodesic flow
is complete in restriction to $A_L$. \item If $M$ is compact then
the geodesic flow is complete and for any leaf $L$ the induced
Riemannian metric $\prs_L$ is complete.\end{enumerate}\end{co}

We will now construct an analogous of the Sasaki metric on $A$ and
study the divergence of the geodesic flow with respect to this
metric. Actually, the Sasaki metric is not defined on $A$ but only on $A_L$ where $L$ is a leaf of the characteristic foliation.\\

 Let $p:A\too M$ be a Riemannian Lie algebroid with anchor
$\#$. Fix a leaf $L$, consider $p_L:A_L\too L$ and put $\V A_L={\mathrm Ker}dp_L$.

For any $a\in A_L$,  we consider  the
subspace $\h^\perp A_L$ of $T_a A_L$ consisting of the tangent vectors $V_a$ such
that there exists an horizontal $A$-path $\al:[0,1]\too \G_L^\perp$
satisfying $p(\al(0))=p(a)$ and
$$V_a=\frac{d}{dt}_{|t=0}\tau_\al^t(a),$$where  $\tau_\al$ is the parallel transport along $\al$. We have \begin{equation} TA_L=\V A_L\oplus \h^\perp A_L.\end{equation}
Indeed, we
 define $$K:TA_L\too A_L$$ as follows. Fix $a\in A_L$ and $Z\in T_a
A_L$ and choose $\be:[0,1]\too A_L$ such that $\be(0)=a$ and
$\dot\be(0)=Z$. There exists an unique horizontal $A$-path
$\al:[0,1]\too \G_L^\perp$ with  the base path $p\circ\be(t)$. Put
$$K(Z)=(\D^\al\be)(0).$$

It is easy to check that $K$ is well-defined, ${\mathrm
Ker}K=\h^\perp A_L$, for any $Z\in\V A_L$ $K(Z)=Z$ and deduce $(17)$.

Let $(x_1,\ldots,x_l)$ be a system of local coordinates on an open set $U$ in $L$
and $(a_1,\ldots,a_r)$ is a basis of local sections (over $U$) of $A_L$. This
defines a system of coordinates  $(x_1,\ldots,x_l,\mu_1,\ldots,\mu_r)$
on $A_L$ and  if
$$Z=\sum_{j}b_j\partial_{x_j}+\sum_{j}Z^j\partial_{\mu_j}$$then
\begin{equation}K(Z)=\sum_{l}\left(Z^l+\sum_{i,j}\al_i\mu_j\Ga_{ij}^l\right)a_l\end{equation}where
$dp_L(Z)=\#(\sum_{i}\al_ia_i)$ and $\sum_{i}\al_ia_i\in\G_L^\perp$.

\begin{rem} In general, the geodesic vector field  does not lies in ${\mathrm Ker} K$. Indeed, one can check easily
that for any $a\in A_L$
\begin{equation}K(X_E(a))=-\D_{a^v}a.\end{equation}
Moreover, for any $a\in A$, the geodesic $\phi_t(a)$ splits
$\phi_t(a)=\phi_t^v(a)+\phi_t^h(a)$ and the path
$\al(t)=\phi_t^h(a)$ is an $A$-path and the vector field
$$X_E^h(a)=\frac{d}{dt}_{|t=0}\tau_\al^t(a)$$is the horizontal
component of $X_E$.\end{rem}

We define the Sasaki metric on $A_L$ by
$$g_L(Z_a,Z_a)=<d_a p(Z_a),d_a p(Z_a)>_L+<K(Z_a),K(Z_a)>.$$

The projection $p_L:A_L\too L$ becomes a Riemannian submersion. We consider now
  the Liouville vector field $\overrightarrow{r}$ on $A_L$ which is the vector field generating the flow $\phi_t(a)=e^ta$.
By direct computation one can get
\begin{equation}
[\overrightarrow{r},X_E]=X_E.\end{equation} From this relation,
one deduce that $X_E$ preserves the Riemannian volume on $A_L$
associated to $g_L$ if and only if $X_E$ preserves the Riemannian
volume of the restriction of $g_L$ to the spheres bundle
$UA_L=\{a\in A_L;<a,a>=1\}$. Let us compute the divergence of the geodesic vector field with respect to $g_L$.

\begin{th}The divergence the geodesic vector field $X_E$ with respect to the Sasaki metric
$g_L$ is given by
\begin{equation} div(X_E)(a)={\mathrm Tr}ad_{a^v}+<a^h,N>\end{equation}where
$ad_{a^v}:\G_{p(a)}\too\G_{p(a)}$, $b\too[a^v,b]$ and
$N=\sum_{i}T_{b_i}b_i$ where $(b_1,\ldots,b_{s})$ is any
orthonormal basis of $\G_{p(a)}$ and $T$ is the O'Neill tensor defined in 3.1.
\end{th}

{\bf Proof.}  Denote by $l$ the dimension of $L$ and choose a system of  local coordinates $(x_1,\ldots,x_l)$ in some open set $U$ of $L$.
Choose  $(a_1,\ldots,a_l)$ an othonormal basis of  sections of $\G_L^\perp\too U$
and $(b_1,\ldots,b_{r-l})$ an orthonormal basis of sections of $\G_L\too U$. We get a system of coordinates $(x,\mu)$ in $A_L$. Put, for
any $i=1,\ldots,l$,
$$\#(a_i)=\sum_{j}p^{ij}\partial_{x_j}\quad\mbox{and}\quad
Z^i=\sum_{j}p^{ij}\partial_{x_j}-
\sum_{l}(\sum_{j}\mu_j\Ga_{ij}^l)\partial_{\mu_l}.$$ By using $(18)$, one can check easily that
$K(Z^i)=0$ and $K(\partial_{\mu_i})=a_i$ for $i=1,\ldots,l$ and
$K(\partial_{\mu_i})=b_i$ for $i=l+1,\ldots,l-r$. Moreover $(Z^1,\ldots,Z^l,\partial_{\mu_1},\ldots,\partial_{\mu_r})$ is an
orthonormal frame of $g_L$ and hence
$$div(X_E)=\sum_{i}g_L([Z^i,X_E],Z^i)+\sum_{j}g_L([\partial_{\mu_j},X_E],\partial_{\mu_j}).$$
Recall that
$$X_E=\sum_{i,k=1}^lp^{ki}\mu_k\partial_{x_i}-
\sum_{j,s,t}\mu_s\mu_t\Ga_{st}^j\partial_{\mu_j}.$$

 So, for $1\leq j\leq l$,
\begin{eqnarray*}
[\partial_{\mu_j},X_E]&=&\sum_{i}p^{ji}\partial_{x_i} - \sum_{i,t}\mu_t
(\Ga_{jt}^i+\Ga_{tj}^i)\partial_{\mu_i},\\
g_L([\partial_{\mu_j},X_E],\partial_{\mu_j})&=&<K([\partial_{\mu_j},X_E]),K(\partial_{\mu_j})>
\stackrel{(18)}=-\sum_{t}\mu_t\Ga_{tj}^j=0,
\end{eqnarray*}since $\Ga_{tj}^j=<\D_{a_t}a_j,a_j>=-<a_j,\D_{a_t}a_j>.$

For $j\geq l+1$
\begin{eqnarray*}
[\partial_{\mu_j},X_E]
&=&-\sum_{i,t}\mu_t
(\Ga_{jt}^i+\Ga_{tj}^i)\partial_{\mu_i},\\
g_L([\partial_{\mu_j},X_E],\partial_{\mu_j})&=&-
\sum_{t}\mu_t\Ga_{jt}^j.\\
\end{eqnarray*}
Hence
\begin{eqnarray*}
\sum_{j}g_L([\partial_{\mu_j},X_E],\partial_{\mu_j})&=&-\sum_{j\geq
l+1}\sum_{t}\mu_t\Ga_{jt}^j\\
&=&-\sum_{j\geq
l+1}\sum_{t=1}^l\mu_t<\D_{b_j}a_t,b_j>-\sum_{j\geq
l+1}\sum_{t\geq l+1}\mu_t<\D_{b_j}b_t,b_j>\\
&=&<a^h,\sum_{j\geq l+1}\D_{b_j}b_j>-\sum_{j\geq
l+1}<\D_{b_j}a^v,b_j>\\
&=&<a^h,\sum_{j\geq l+1}T_{b_j}b_j>-\sum_{j\geq
l+1}<[b_j,a^v],b_j>\\
&=&<a^h,N>+{\mathrm Tr}ad_{a^v}.\end{eqnarray*}

On the other hand, one can see easily
that\begin{eqnarray*}X_E&=&\sum_{k=1}^l\mu_kZ^k-\sum_{j=1}^r
\sum_{s\geq
l+1,t}\mu_s\mu_t\Ga_{st}^j\partial_{\mu_j}=\sum_{k=1}^l\mu_kZ^k+V.\end{eqnarray*} Note that $V$ is vertical and since, for any $i=1,\ldots,l$, $Z^i$ is basic (with respect to the Riemannian submersion $p_L:A_L\too L$) then $[Z^i,V]$ is vertical. Note also that, for any $i,k=1,\ldots,l$, $dp_L([Z^i,Z^k])=\#([a_i,a_k])$. Hence
\begin{eqnarray*}
\sum_{i}g_L([Z^i,X_E],Z^i)&=&\sum_{i,k}g_L([Z^i,\mu_kZ^k],Z^i)\\
&=&\sum_{i}Z^i(\mu_i)+\sum_{i,k}\mu_k<\#([a_i,a_k]),\#a_i>_L\\
&=&\sum_{i}Z^i(\mu_i)+\sum_{i,k}\mu_k<[a_i,a_k]^h,a_i>\\
&=&\sum_{i}Z^i(\mu_i)+\sum_{i,k}\mu_k<[a_i,a_k],a_i>\\
&=&-\sum_{i,k}\mu_k\Ga_{ik}^i+\sum_{i,k}\mu_k\Ga_{ik}^i=0.
\end{eqnarray*}

Finally, we get
the desired formula.$\Box$

\begin{co} The geodesic flow preserves the Riemannian volume on $A_L$ if
and only if $N=0$ and $\G_x$ is unimodular for some $x\in
L$.\end{co}

\begin{rem}\begin{enumerate}\item If $A=TM$ then $\G=\{0\}$ and $div(X_E)=0$ and hence we recover the classical theorem of Liouville.
\item If $A$ is a Lie algebra $\G$ endowed with a scalar product $\prs$. The geodesic vector field $X_E$ is given, in any system of linear coordinates $(x_1,\ldots,x_n)$, by
    $$X_E=-\sum_{i,s,t}x_sx_t\Ga_{st}^i\partial_{x_i}$$and the Sasaki metric is the flat Riemannian metric $\prs$. From Theorem 3.2 we deduce that $div X_E(a)={\mathrm Tr}ad_a$. Hence $X_E$ is divergence free if and only if $\G$ is unimodular.

\end{enumerate}\end{rem}

We will now establish the first and the second variation formulas
in the context of Riemannian Lie algebroids.

Let $p:A\too M$ be a Riemannian Lie algebroid with anchor $\#$.
For any $A$-path $\al:[0,1]\too A$ we define the energy and the
length of $\al$, respectively, by
$${\mathbf E}(\al)=\frac12\int_0^1<\al(t),\al(t)>dt\quad\mbox{and}\quad
\Li(\al)=\int_0^1\sqrt{<\al(t),\al(t)>}dt.$$

For any $m,q$ lying in the same leaf of the characteristic
foliation, we denote by $\Om_{mq}$ the set of $A$-path $\al$ such
that $p(\al(0))=m$ and $p(\al(1))=q$.
\begin{pr}{\bf (First variation formulas)} Let $p:A\too M$ be a
Riemannian Lie algebroid. Then:
\begin{enumerate}\item For any variation of $A$-paths  $\al:[0,1]\times [0,1]\too A$
and for any $\be$ a transverse variation to $\al$, one has
\begin{eqnarray*}
 \frac{d}{d\e}{\mathbf
 E}(\al)&=&<\be(\e,1),\al(\e,1)>-<\be(\e,0),\al(\e,0)>-\int_0^1< \be,\D_t
\al>dt\\&&- \int_0^1<\De(\al,\be),\al>dt.\end{eqnarray*} \item The
$h$-critical points of ${\mathbf E}:\Om_{mq}\too\reel$, namely the
$A$-paths $\al_0$ such that $$\frac{d}{d\e}{\mathbf
 E}(\al)_{|\e=0}=0$$ for any $A$-homotopy $\al$ in  $\Om_{mq}$ starting at $\al_0$,
 are
 geodesics.
 \item  For any
 variation of $A$-paths $\al$ such that $\al_0$ is parameterized with
 arc-length,
  $$\frac{d}{d\e}{\mathbf
 E}(\al)_{|\e=0}=\frac{d}{d\e}{\Li
 }(\al)_{|\e=0}.$$
 \item An $A$-path $\al_0\in\Om_{mq}$ is $h$-critical for $\Li$, namely
 $$\frac{d}{d\e}{\Li}(\al)_{|\e=0}=0$$ for any $A$-homotopy in
 $\Om_{mq}$  starting at $\al_0$, if and only if there exists
 a change of parameter $\mu$ such that the $A$-path $\wi
 \al_0=\mu'\al_0(\mu)$ is a geodesic.

\end{enumerate}\end{pr}

{\bf Proof.} \begin{enumerate}\item Let  us compute
$\frac{d}{d\e}{\mathbf E}(\al)$. We have
\begin{eqnarray*}
\frac{d}{d\e}{\mathbf
E}(\al_\e)&=&\frac12\frac{d}{d\e}\int_0^1<\al,\al>dt
=\frac12\int_0^1\frac{d}{d\e}<\al,\al>dt
=\int_0^1<\D_\e \al,\al>dt\\
&=&\int_0^1<\D_t \be,\al>dt
-\int_0^1<\De(\al,\be),\al>dt\quad(\mbox{Proposition 2.2})\\
&=&\int_0^1\partial_t(< \be,\al>)dt -\int_0^1(< \be,\D_t \al>)dt-
\int_0^1<\De(\al,\be),\al>dt\\
&=&< \be(\e,1),\al(\e,1)>-< \be(\e,0),\al(\e,0)> -\int_0^1<
\be,\D_t \al>dt\\&&-
\int_0^1<\De(\al,\be),\al>dt.\\
\end{eqnarray*}
Analogously one can get
\begin{eqnarray}
\frac{d}{d\e}\Li(\al)=\int_0^1|\al|^{-1/2}\partial_t(< \be,\al>)dt
-\int_0^1|\al|^{-1/2}(< \be,\D_t \al>)dt\nonumber\\-
\int_0^1|\al|^{-1/2}<\De(\al,\be),\al>dt.\end{eqnarray}

\item Let $\al_0$ be geodesic and let $\al$ be an $A$-homotopy with
fixed end-point starting at $\al_0$. Then there exists a
transverse variation $\be$ to $\al$ such that
$\be(\e,0)=\be(\e,1)=0$ and $\De(\al,\be)=0$. Hence from 1., we
get
$$\frac{d}{d\e}_{|\e=0}{\mathbf E}(\al)=0.$$
Conversely, suppose that $\al_0$ is an $A$-path which is a
$h$-critical point of ${\mathbf E}:\Om_{mq}\too\reel$. Consider
the $\al_0$-section $\be_0(t)=f(t)\D_t\al_0$ where
$f:[0,1]\too\reel$ is a smooth function such that $f(0)=f(1)=0$.
According to Lemma 2.1, there exists an $A$-homotopy $\al$ with fixed
end-points and starting at $\al_0$ and such the corresponding
transverse variation $\be$ satisfies $\be(0,t)=\be_0(t)$. By
applying the formula in 1., we get
$$0=\int_0^1f(t)<\D_t\al_0,\D_t\al_0>dt$$ and hence $\D_t\al_0=0$
which means that $\al_0$ is a geodesic. \item This is a
consequence of $(22)$ and $|\al_0|=1$. \item Immediate from 2. and
3. $\Box$\end{enumerate}

\begin{pr}{\bf (Second variation formulas)} Let $p:A\too M$ be a
Riemannian Lie algebroid. Then the following assertions hold.
\begin{enumerate}\item For any  variation of $A$-paths $\al$ such that $\al_0$ is a geodesic
and for any $\be$  a transverse variation to $\al$ such that
$\De(\al,\be)=0$, one has
\begin{eqnarray*}
 \frac{d^2}{d\e^2}{\mathbf
 E}(\al)_{|\e=0}&=&<\D_\e \be(0,1),\al(0,1)>-<\D_\e \be(0,0),\al(0,0)>
 \\&&+\int_0^1<  \D_t \be_0,\D_t\be_0 >dt
 +\int_0^1<\be_0,R(\al_0,\be_0)\al_0>dt.\end{eqnarray*}
\item Let $\al$ be an $A$-homotopy of $A$-paths such that $\al_0$ is a
geodesic and let $\be$ be the corresponding transverse variation.
One has
\begin{eqnarray*}
 \frac{d^2}{d\e^2}{\mathbf
 E}(\al)_{|\e=0}&=&
 \int_0^1<  \D_t \be_0,\D_t\be_0
 >dt+\int_0^1<\be_0,R(\al_0,\be_0)\al_0>dt.\end{eqnarray*}
\item Let $\al$ be a variation of $A$-paths such that $\al_0$ is a
geodesic parameterized by arc length and let $\be$ a transverse
variation to $\al$ such that $\De(\al,\be)=0$. One has
\begin{eqnarray*}
 \frac{d^2}{d\e^2}{\Li
 }(\al)_{|\e=0}&=&<\D_\e \be(0,1),a(0,1)>-<\D_\e \be(0,0),\al(0,0)>
 \\&&+\int_0^1<  \D_t \be_0,\D_t\be_0
 >dt+\int_0^1<\be_0,R(\al_0,\be_0)\al_0>dt\\&&-\int_0^1<\al_0,\D_t \be_0>dt.\end{eqnarray*}
\item Let $\al$ be an $A$-homotopy of $A$-paths such that $\al_0$ is a
geodesic parameterized by arc length and let $\be$ be the
corresponding transverse variation. One has
\begin{eqnarray*}
 \frac{d^2}{d\e^2}{\Li }(\al)_{|\e=0}&=&
 \int_0^1<  \D_t \be_0,\D_t\be_0
 >dt+\int_0^1<\be_0,R(\al_0,\be_0)\al_0>dt\\&&-\int_0^1<\al_0,\D_t \be_0>dt.\end{eqnarray*}

\end{enumerate}\end{pr}

{\bf Proof.}\begin{enumerate}\item
 From the first variation formula, we have
\begin{eqnarray*}
 \frac{d}{d\e}{\mathbf
 E}(\al)&=&<\be(\e,1),\al(\e,1)>-<\be(\e,0),\al(\e,0)>-\int_0^1< \be,\D_t
\al>dt.\\\end{eqnarray*}

Then
\begin{eqnarray*}
 \frac{d^2}{d\e^2}{\mathbf
 E}(\al)&=&<\D_\e \be(\e,1),\al(\e,1)>+< \be(\e,1),\D_\e \al(\e,1)>
 \\&&-<\D_\e \be(\e,0),\al(\e,0)>-<\be(\e,0),\D_\e \al(\e,0)>\\&&
 -\int_0^1< \D_\e \be,\D_t\al>dt-\int_0^1<  \be,\D_\e\D_t\al>dt.\\
 \int_0^1<  \be,\D_\e\D_t\al>dt&=&\int_0^1<  \be,\D_t\D_\e
 \al>dt+\int_0^1<\be,R(\be,\al)\al>dt\\
 &=&\int_0^1\partial_t(<  \be,\D_\e
 \al>)dt-\int_0^1<  \D_t \be,\D_\e
 \al>dt\\&&+\int_0^1<\be,R(\be,\al)\al>dt\\
 &=&<  \be(\e,1),\D_\e
 \al(\e,1)>-<  \be(\e,0),\D_\e
 \al(\e,0)>\\
 &&-\int_0^1<  \D_t \be,\D_t\be
 >dt+\int_0^1<\be,R(\be,\al)\al>dt\\
 \end{eqnarray*}

 Hence
 \begin{eqnarray*}
 \frac{d^2}{d\e^2}{\mathbf
 E}(\al)&=&<\D_\e \be(\e,1),\al(\e,1)>-<\D_\e \be(\e,0),\al(\e,0)>-\int_0^1< \D_\e
 \be,\D_t\al>dt\\
 &&+\int_0^1<  \D_t \be,\D_t\be
 >dt+\int_0^1<\be,R(\al,\be)\al>dt.\end{eqnarray*}
\item In this situation, we have $\D_\e \be(\e,1)=\D_\e
\be(\e,0)=\D_t\al=0$ and the formula follows.

3. and 4. are left to the reader.$\Box$\end{enumerate}

As an application of  Proposition 3.3, we give now a description of the geodesics of a left invariant Riemannian metric on a Lie group using the geodesics of its Lie algebra considered as a Riemannian Lie algebroid.\\

Let $G$ be a Lie group and $\G=T_eG$ its Lie algebra. For any $u\in\G$, we denote by $u^+$ the associated left invariant vector field on $G$. Suppose that $G$ is endowed with a left invariant Riemannian metric $g$ and put $\prs=g_e$. If we think to $\G$ as a Lie algebroid, $(\G,\prs)$ is a Riemannian Lie algebroid and we will explain how one can construct the geodesics of $(G,g)$ from the geodesics of $(\G,\prs)$. Choose a basis $(e_1,\ldots,e_n)$  of $\G$ and put $g_{ij}=<e_i,e_j>$. Recall that the geodesics of $(\G,\prs)$ are the integral curves of the geodesic vector field $X_E$ given in the linear coordinates $(x_1,\ldots,x_n)$ associated to $(e_1,\ldots,e_n)$ by
$$X_E=-\sum_{s,t,j}x_sx_t\Ga_{st}^j\partial_{x_j},$$where $\Ga_{st}^j$ are given by
$$\Ga_{st}^j=\frac12\sum_{l,u}g^{lj}\left(g_{ul}C_{st}^u+g_{ut}C_{ls}^u+g_{us}C_{lt}^u\right).$$ Here $(g^{ij})$ is the inverse matrix of $(g_{ij})$ and $C_{ij}^k$ are given by $[e_i,e_j]=\sum_{u}C_{ij}^ue_u.$
\begin{pr} Let $h\in G$ and $v\in T_hG$. Then the geodesic $\ga:\reel\too G$ of $(G,g)$ satisfying $\ga(0)=h$ and $\dot\ga(0)=v$ is the integral curve passing through $h$ of the time-depending family of left invariant vector fields $(\al^+(t))_{t\in\reel}$ where $\al:\reel\too\G$ is the geodesic of $(\G,\prs)$ satisfying $\al(0)=\left(L_{h^{-1}}\right)_*(v).$\end{pr}

{\bf Proof.} Note first that by invariance the integral curves of $(\al^+(t))_{t\in\reel}$ are complete. Note also that both $(G,g)$ and $(\G,\prs)$ are geodesically complete. Let $\ga:\reel\too G$ be the integral curve of $(\al^+(t))_{t\in\reel}$ satisfying $\ga(0)=h$. We have
$$\dot\ga(0)=\al^+(0)=\left(L_{h}\right)_*(\al(0))=\left(L_{h}\circ L_{h^{-1}}\right)_*(v)=v.$$
We will show that for any $t_1,t_2\in\reel$, the restriction of $\ga$ to $[t_1,t_2]$ is a critical point of the  energy functional ${\mathbf E}_g:\Om\too\reel$ where $\Om$ is the space of smooth curves $\mu:[t_1,t_2]\too G$ such that $\mu(t_1)=\ga(t_1)$ and $\mu(t_2)=\ga(t_2)$.

Let $\wi\ga:[0,1]\times[t_1,t_2]\too G$ be an homotopy with end-fixed points such that $\wi\ga(0,.)=\ga$. It is well-known (see [7]) that the variation $\wi\al:[0,1]\times[t_1,t_2]\too \G$ given by
$$\wi\al(\e,t)=\left(L_{\wi\ga(\e,t)^{-1}}\right)_*\left(\frac{\partial\wi\ga}{\partial t}(\e,t)\right)$$ is a $\G$-homotopy. Moreover, $\wi\al(0,.)=\al$ and, by invariance, ${\mathbf E}_{g}(\wi\ga)={\mathbf E}_{\prs}(\wi\al)$. By applying Proposition 3.3, we get
$\frac{d} {d\e}{\mathbf E}_{\prs}(\wi\al)_{|\e=0}=0$. Thus, $\frac{d }{d\e}{\mathbf E}_g(\wi\ga)_{|\e=0}=0$ and, by applying the classical result on geodesics of Riemannian metric we deduce that $\ga$ is a geodesic.$\Box$

\begin{rem} If the Riemannian metric $g$ is bi-invariant then $\Ga_{ij}^k=\frac12 C_{ij}^k$ and hence $X_E$ vanishes identically. We deduce from Proposition 3.5 that the geodesic of $(G,g)$ passing through $h\in G$ and with initial velocity $v\in T_hG$ is the integral curve (passing through $h$) of the left invariant vector field $(L_{h^{-1}*}(v))^+$.\end{rem}

Let us define now Jacobi sections along a geodesic.

\begin{Def}Let $A$ be a Riemannian Lie algebroid and $\al:[0,1]\too A$ a geodesic.
 A Jacobi $\al$-section  is an $\al$-section $\be$
which satisfies
$$\be''-R(\al,\be)\al=0,$$where $\be'$ is the derivative of $\be$ along $\al$ and so on.
\end{Def}
\begin{pr} Let $\al:[0,1]\too A$ be a geodesic in a Riemannian Lie algebroid
$A$. Then
 for any $a,b\in A_{p(\al(0))}$ there exists
one and only one Jacobi $\al$-section  such that $\be(0)=a$
and $\be'(0)=b$. If $\be(0)=0$ and $\be'(0)=k\al(0)$ then
$\be(t)=kt\al(t)$ for any $t$. If $\be(0)$ and $\be'(0)$ are
orthogonal to $\al(0)$, then $\be(t)$ is orthogonal to $\al(t)$ for
any $t$. In particular the vector space of Jacobi $\al$-sections has
dimension $2r$ and the subspace of Jacobi $\al$-sections which are
normal to $\al$ has dimension $2(r-1)$.
\end{pr}

{\bf Proof.} Take an orthonormal basis $(a_1,\ldots,a_r)$ of
$A_{p(\al(0))}$ such that $a_1=k\al(0)$. The parallel transport along
$\al$ of the vectors $a_i$ gives a basis of orthonormal $\al$-sections
$(s_1,\ldots,s_r)$ with $s_1=k\al$.
Every Jacobi $\al$-section $\be$  is a linear combination of
$s_i$, say $\be=\sum_{i}y_is_i$, whose coefficients satisfy the
differential system
$$y_i''-\sum_{j=2}^r<R(\al,s_j)\al,s_i>y_j=0.$$For given initial conditions
$\be(0)=a$ and $\be'(0)=b$, the existence and uniqueness of $\be$
come from standard results about linear differential systems.

If $\be(0)=0$ and $\be'(0)=k\al(0)$ then $\be(t)=kt\al(t)$ since
$\be''(t)=0$.

The condition $\be(0)$ and $\be'(0)$ to be orthogonal to $\al$
means that $y_1(0)=0$ and $y'_1(0)=0$. In that case $y_1(t)=0$ for
any $t$, since $y''(t)=0$.$\Box$
\begin{pr} Let $\al_0:[0,1]\too A$ be a geodesic, and $\al$ be a
variation of $\al_0$ such that all $A$-paths $\al(\e,.)$ are
geodesics. Then, for any transverse variation $\be$ of $\al$ such
that $\De(\al,\be)=0$,  $\be_0$ is a Jacobi $\al_0$-section. Conversely, every Jacobi $\al_0$-section can be obtained in this
way.\end{pr}

{\bf Proof.} We have
$$\be''_0(t)=\D_t\D_t\be(0,t).$$Performing the two exchanges of $t$ and
$\e$, we get from Proposition 2.2
$$\be''_0(t)=\D_t\D_\e \al(0,t)=\D_\e\D_t\al(0,t)+R(\al_0,\be_0)\al_0.$$ Since the
$A$-paths $\al_\e$ are geodesics, the first term vanishes and we
get
$$\be''_0=R(\al_0,\be_0)\al_0.$$

Conversely, take a Jacobi $\al_0$-section $b$  and the
geodesic $c$ such that $c(0)=b(0)$. Take parallel sections
$s_0$ and $s_1$ along $c$ such that $s_0(0)=\al_0(0)$ and
$s_1(0)=b'(0).$ Set
$$s(\e)=s_0(\e)+\e s_1(\e)\quad\mbox{and}\quad
\al(\e,t)=\phi_t(s(\e)),$$where $\phi_t$ is the geodesic flow.
Consider the transverse variation $\be$ to $\al$ such that
$\be(\e,0)=c(\e)$ and $\De(\al,\be)=0$. We will show that
$\be(0,.)$ and $b$ coincide. Remark first that these two
$\al_0$-sections  satisfy the same differential equation
namely
$$y''-R(\al_0,y)\al_0=0.$$ Since $b(0)=\be(0,0)=c(0)$, let us show that
$\D_t\be(0,0)=b'(0)$. Since $\D_t\be=\D_\e \al$, we have
$\D_\e \al(0,0)$ is the value at 0 of the derivative of the curve
$\al(\e,0)$ along the $A$-path $\be(\e,0)$. Or $\al(\e,0)=s(\e)$
and $\be(\e,0)=c(\e)$ and we get $\D_\e
\al(0,0)=s_1(0)=b'(0)$.$\Box$\\

As the classical case, the Jacobi sections can be used to compute the derivative of the exponential which can be defined as follows. Let $p:A\too M$ be a Riemannian Lie algebroid.
 Fix a point $m\in M$ and denote
by $L$ the leaf containing $m$. We define the exponential
$$exp_m:\U\subset A_m\too L$$where $\U_m=\{\al\in A_m,
\phi_1(\al)\quad\mbox{is defined}\}$ and
$exp_m(\al)=p\circ\phi_1(\al)$ ($\phi$ is the geodesic flow).

\begin{pr} We have
$$d_{a}exp_m(u)=\#(\be(1))$$where $\be$ is the Jacobi section along
$t\mapsto\phi_t(a)$ with initial condition $\be(0)=0$ and
$\be'(0)=u$.\end{pr}

{\bf Proof.} We have
$$d_{a}exp_m(u)=\frac{d}{d\e}_{|\e=0}p(\phi_1(a+\e u)).$$
We consider the variation of geodesics $\al(\e,t)=\phi_t(a+\e
u)$ with fixed initial point. We consider the transverse variation
$\be$ such that $\be(\e,0)=0$ and $\De(\al,\be)=0$. We have that
$\be_0$ is a Jacobi $\al_0$-section such that $\be_0(0)=0$ and
$\#(\be_0(1))=\frac{d}{d\e}_{|\e=0}p(\phi_1(a+\e u))$ by
construction.$\Box$\\

As the classical case, we define the sectional curvature of two linearly independent  vectors $a,b\in A_m$ by
$$K(a,b)=-\frac{<R(a,b)a,b>}{<a,a><b,b>-<a,b>^2}.$$
\begin{pr} Let $p:A\too M$ be  a Riemannian Lie algebroid. If the sectional curvature is everywhere nonpositive then $exp_m$ is a
submersion for every $m\in M$.\end{pr}

{\bf Proof.} Fix $a\in A_m$ and let $\J_0^\al$ be the space of
Jacobi sections $\be$ along $\al(t)=\phi_t(a)$ such that
$\be(0)=0$ ($\phi$ is the geodesic flow). We define the linear application
$$\xi:\J_0^\al\too A_{p(\phi_1(\al_0))}$$ by
$\xi(\be)=\be(1).$
We will show that $\xi$ is injective and hence an isomorphism
since $\dim\J_0^\al=\dim A_{p(\phi_1(\al_0))}$. Suppose that
$\be\in \J_0^\al$ satisfies $\be(1)=0$. The function
$f:[0,1]\too\reel$ given by $f(t)=<\be(t),\be(t)>$ satisfies
\begin{eqnarray*}
f'(t)&=&2<\be'(t),\be(t)>,\\
f''(t)&=&2<\be'(t),\be'(t)>+2<\be''(t),\be(t)>\\
&=&2<\be'(t),\be'(t)>+2<R(\al(t),\be(t))\al(t),\be(t)>.\end{eqnarray*}
Hence $f''\geq0$ and since $f(0)=f(1)=0$ we deduce that $f$
vanishes identically and then $\be=0$. This shows that $\xi$ is
injective and hence an isomorphism. From Proposition 3.8, one can
identify ${\mathrm Ker}d_{a}exp_m$ with
$\xi^{-1}(\G_{p(\phi_1(a))})$ and the proposition
follows.$\Box$

\section{O'Neill's formulas for curvature}

Let $p:A\too M$ be a Riemannian Lie algebroid. The different
curvatures (sectional curvature, Ricci curvature and scalar
curvature) can be defined as the classical case (when $A=TM$). For any leaf $L$,
the short exact sequence
$$0\too\G_L\too A_L\too TL$$ is formally identical to a Riemannian
submersion and hence all formulas on curvature given by O'Neill
are valid in this context. We denote by $K$, $\hat K$ and $\wi K$ respectively, the sectional curvature of the Riemannian metrics $\prs$, the restriction of $\prs$ to $\G_L$ and the induced metric on $L$. The following proposition is a reformulation of Corollary 9.29 pp. 241 in [1].

\begin{pr} Let $\al,\be,s_1,s_2\in\Ga(A_L)$ such that $\al,\be$ are vertical, $s_1,s_2$ are horizontal and $|\al\wedge\be|=1$, $|s_1|=|\al|=1$, $|s_1\wedge
s_2|=1$. Then
\begin{eqnarray*} K(\al,\be)&=&\hat
K(\al,\be)+|T_\al\be|^2-<T_\al\al,T_\be\be>,\\
K(s_1,\al)&=&<(\D_{s_1}T)_\al\al,s_1>-|T_\al
s_1|+|H_{s_1}\al|^2,\\
K(s_1,s_2)&=&\wi K(s_1,s_2)-3|H_{s_1}s_2|^2.\end{eqnarray*}\end{pr}

The last formula says that the leaves carry "more curvature" than the Lie algebroid and by applying Mayer theorem we get:
\begin{pr} Let $A\too M$ be a complete Riemannian algebroid and let $L$ be
a leaf of the characteristic foliation such that for any linearly independent  horizontal sections $s_1,s_2$ over
$L$, $K(s_1,s_2)\geq k$. Then
 $diam L\leq\frac{\pi}{\sqrt{k}}$ and hence $L$ is compact.\end{pr}

 There is another case when one can apply Mayer theorem. Consider a Riemannian Lie algebroid $p:A\too M$ such that the O'Neill tensor $T$ vanishes and fix a leaf $L$ and denote by $r$ and $\wi r$ respectively the Ricci curvature of the Riemannian metrics $\prs$ and $\prs_L$. The formula $9.36c$ pp.244 in $[1]$ applies in our context and gives
 $$r(s_1,s_2)=\wi r(\#(s_1),\#(s_2))-2\sum_{i=1}^l<H_{s_1}a_i,H_{s_2}a_i>$$where $(a_1,\ldots,a_l)$ is any orthonormal basis of $\G_L^\perp$. By applying Mayer theorem we get:
 \begin{pr} Let $A\too M$ be a complete Riemannian algebroid such that $T=0$ and let $L$ be
a leaf of the characteristic foliation such that there exists a constant $k$ such that the restriction of $r$ to $\G_L^\perp$ satisfies $$r\geq(n-1)k^{-2}\prs.$$
 Then $diam L\leq\frac{\pi}{\sqrt{k}}$ and hence $L$ is compact.\end{pr}

\section{Integrability of Riemannian Lie algebroids}
In this section, we study the integrability of Riemannian Lie algebroids. We show that a Riemannian Lie algebroid such that the O'Neill tensor $H$ vanishes is integrable. We show also that a complete Riemannian Lie algebroid with nonpositive sectional curvature is intergrable and it is diffeomorphic to its Weinstein Lie groupoid. This result is a generalization of Hadamard-Cartan theorem.\\

A groupoid is a small category $\C$ in which all the arrows are
invertible. We shall write $M$ for the set of objects of $\C$,
while the set of arrows of $\C$ will be denoted by $\C$. We shall
often identify $M$ with the subset of units of $\C$. The structure
maps of $\C$ will be denoted as follows: ${\mathbf s},{\mathbf
t}:\C\too M$ will stand for the source map, respectively the
target map, $m:\C^2=\{(g,h);{\mathbf s}(g)={\mathbf t}(h)\}\too\C$
the multiplication map $(m(g,h)=gh)$, $i:\C\too\C_1$
$(i(g)=g^{-1})$ for the inverse map and $u:M\too \C$ $(u(x)=1_x)$
for the unit map. Given $g\in\C$, the right multiplication by $g$
is only defined on the ${\mathbf s}$-fiber at ${\mathbf t}(g)$,
and induces a bijection
$$R_g:{\mathbf s}^{-1}({\mathbf t}(g))\too{\mathbf s}^{-1}({\mathbf
s}(g)).$$

A Lie groupoid is a groupoid $\C$, equipped with the structure of
smooth manifold both on the $\C$ and on the $M$ such that all the
structure maps are smooth and ${\mathbf s}$ and ${\mathbf t}$ are
submersions.

The construction of a Lie algebra of a given Lie group extends to
Lie groupoids. Explicitly, if $\C$ is a Lie groupoid, the vector
bundle $T^{{\mathbf s}}\C={\mathrm Ker}(d{\mathbf s})$ over $\C$
of ${\mathbf s}$-vertical tangent vectors pulls back along
$i:M\too\C$ to a vector bundle $A$ over $M$. This vector bundle
has the structure of a Lie algebroid. Its anchor $\#:A\too TM$ is
induced by the differential of the target map, $d{\mathbf
t}:T\C\too TM$. The sections of $A$ over $M$ can be identified by
the space of right invariant ${\mathbf s}$-vertical vector fields
which induce a Lie bracket on the space of sections of $A$. With
this construction in mind, one can see that a Riemannian structure
on $A$ is equivalent to the data of a Riemannian metric on any
${\mathbf s}$-fiber such that, for any $g\in\C$, $R_g:{\mathbf
s}^{-1}({\mathbf t}(g))\too{\mathbf s}^{-1}({\mathbf s}(g))$ is an
isometry. In this case, for any $x\in M$, ${\mathbf t}:{\mathbf
s}^{-1}{(x)}\too L_x$ is a Riemannian submersion where the leaf $L_x$ is
endowed with the metric defined in 3.1.

A Lie algebroid $A$ is called integrable if it is isomorphic to
the Lie algebroid associated to a Lie groupoid. In [7], Crainic
and Fernandes give a final  solution to the problem of  integrability of Lie
algebroids. They show that the obstruction to integrability can be
controlled by two computable quantities.

The following proposition is a direct application of
Crainic-Fernandes results on integrability.
\begin{th} Let $p:A\too M$ be a Riemannian Lie algebroid such that
$H=0$. Then $A$ is integrable.\end{th}

{\bf Proof.} For any leaf $L$, the vanishing of $H$ implies, according to $(11)$, that
the space of sections of  $\G_L^\perp\too L$ is a Lie subalgebra
of $\Ga(A_L)$ and hence there is a splitting $\sigma:TL\too A_L$
of the anchor, which is compatible with the Lie bracket. By
applying Corollary 5.2 in [7], we get the result.$\Box$\\

There is a large class of Lie algebroids for which one can apply this result.
Let $(M,\pi)$ be a Poisson manifold. The cotangent bundle $T^*M$
carries a structure of a Lie algebroid where the anchor is the
contraction by $\pi$, $\pi_\#:T^*M\too TM$ and the Lie bracket is
given by the Koszul bracket
$$[\al,\be]=\Li_{\pi_\#(\al)}\be-\Li_{\pi_\#(\be)}\al-d\pi(\al,\be)$$where
$\al,\be\in\Om^1(M)$. Let $\prs$ be a Riemannian structure in
$T^*M$. In [3], the author studied the triple $(M,\pi,\prs)$ such
that $\pi$ is parallel with respect the  Levi-Civita
( $T^*M$-connection $\D$. A triple $(M,\pi,\prs)$ satisfying $\D\pi=0$ is  called
Riemann-Poisson manifold. The condition $\D\pi=0$ implies that
${\mathrm Ker}\pi_\#$ is invariant by parallel transport and hence
$\D$ is strongly compatible with the Lie algebroid structure of
$T^*M$. By Proposition 3.2 we deduce that $H=0$. So we get the
following result.
\begin{co} Let $(M,\pi,\prs)$ be a Riemann-Poisson manifold. Then
the Lie algebroid structure of $T^*M$ associated to $\pi$ is
integrable.\end{co}

Recall that the Weinstein groupoid of a Lie algebroid $A$ is the
set $\C(A)$ consisting of  $A$-homotopy classes of $A$-paths (see
[7] for detail). The groupoid $\C(A)$ is a topological  groupoid and it carries a
structure of Lie groupoid if and only if $A$ is integrable.

Given a Riemannian Lie algebroid $p:A\too M$, we define the
exponential from $\U A=\{a\in A, \phi_1(a)\quad\mbox{is
defined}\}$ to the Weinstein groupoid $\C(A)$ as follows: $Exp:\U
A\too\C(A)$ maps $a$ to the $A$-homotpoy class of the geodesic
$\phi_t(a):[0,1]\too A$.

A proof of the following theorem will be given in a furthercoming paper.

\begin{th} Let $p:A\too M$ be a complete Riemannian Lie algebroid such
that the sectional  curvature is nonpositive. Then:
\begin{enumerate}\item for any leaf $L$, $exp_L:A_m\too L$ is a
locally trivial fibration,\item $p:A\too M$ is integrable and
$Exp:A\too\C(A)$ is a diffeomorphism.\end{enumerate}\end{th}

 {\bf References}\bigskip

 [1] {\bf A. Besse,} Einstein manifolds, Springer-Verlag,
Berlin-Hiedelberg-New York (1987).

 [2] {\bf
M. Boucetta, } {\it Compatibilit\'e des structures
pseudo-riemanniennes et des structures de Poisson}{ C. R. Acad.
Sci. Paris, {\bf t. 333}, S\'erie I, (2001) 763--768.}

[3] {\bf M. Boucetta, } {\it Poisson manifolds with compatible
pseudo-metric and pseudo-Riemannian Lie algebras}, { Differential
Geometry and its Applications, {\bf Vol. 20, Issue 3} (2004),
279-291.}

[4] {\bf A. Cannas da Silva and A. Weinstein,} {\it Geometric
models for noncommutatve algebras,} Berkeley Mathematics Lecture
Notes 10, Amer. Math. Soc., Providence, 1999.

[5] {\bf P. Cartier,} {\it Groupo\"ides de Lie et leurs
alg\'ebro\"ides}, S\'eminaire Bourbaki 2007-2008, no. 987.

[6] {\bf Cort\`es J., De Le\`on M., Marrero J.C., Martin De Diego and Martinez E.}, {\it A survey of Lagrangian mechanics and control on Lie algebroids and groupoids}, Int. Jour. on Geom. Metho. in Math. Phy., To appear.

 [7] {\bf M. Crainic and  R. Fernandes,  }
  {\it Integrability of Lie brackets,} Ann. of Math. (2) 157 (2003),
  no. 2, 575--620.

  [8] {\bf M. Crainic and R. Fernandes, }
  {\it Integrability of Poisson brackets,} J. Differential Geom. 66 (2004),
  no. 1, 71-137.

[9] {\bf  R. L. Fernandes}, {\it Lie algebroids, holonomy and
characteristic classes}, Adv. in Math. {\bf 170} (2002), 119-179.

  [10] {\bf R. Fernandes, } {\it  Connections in Poisson geometry. I.
  Holonomy and invariants,} J. Differential Geom. {\bf 54} (2000), no. 2, 303--365.

  [11] {\bf K. Grabowska, J. Grabowski and P. Urba\'nski},
  {\it Geometrical Mechanics on algebroids,} Int. J. Geom. Meth.
  Mod.
  Phys. {\bf 3} (2006), 559-575.

[12] {\bf Hawkins E.,} {\it Noncommutative Rigidity}, Commun.
Math. Phys. {\bf 246} (2004) 211-235.

[13] {\bf Hawkins E.,} {\it The structure of noncommutative
deformations,} J. Diff. Geom. {\bf 77}, 385-424 (2007).

 [14] {\bf K. Mackenzie,} {\it Lie Groupoids and
Lie Algebroids in Differential Geometry}, London Math. Soc.
Lecture Notes Ser. {\bf 124}, Cambridge Univ. Press, Cambridge,
1987.

[15] {\bf B. O'Neill,} {\it The fundamental equations of a
submersion}, Mich. Math. J. {\bf 13}, 459-469 (1966).

 [16] {\bf H. Sussmann,} {\it Orbits of families of vector
fields and integrability of distributions}, Trans. Amer. Math.
Soc. {\bf 180} (1973), 171-188.

[17] {\bf Vaisman I.}, {\it Lecture on the geometry of Poisson
manifolds}, Progr. In Math. {\bf Vol. 118}, Birkhausser, Berlin,
(1994).

[18] {\bf A. Weinstein,} {\it Lagrangian mechanics and groupoids,}
Fields Inst. Comm. {\bf 7} (1996), 207-231.\bigskip

Mohamed Boucetta\\
Facult\'e des Sciences et Techniques \\
BP 549 Marrakech, Morocco.
\\
Email: {\it mboucetta2@yahoo.fr }

\end{document}